\documentclass[12pt]{amsart}
\usepackage{amssymb}

\newtheorem{thm}{Theorem}[section]
\newtheorem{prob}[thm]{Problem}

\newtheorem{conj}[thm]{Conjecture}
\newtheorem{issue}{Issue}

\theoremstyle{definition}

\theoremstyle{remark}

\newcommand{\x}{\times}

\newcommand{\gpbl}{{\mbox{\textit{\tiny gp}}}}

\renewcommand{\b}{{\mathfrak b}}

\renewcommand{\d}{{\mathfrak d}}
\newcommand{\p}{{\mathfrak p}}

\newcommand{\g}{\mathfrak{g}}
\newcommand{\ADD}{{\mathsf   {ADD}}}
\newcommand{\COV}{{\mathsf   {COV}}}
\newcommand{\NON}{{\mathsf   {NON}}}
\newcommand{\COF}{{\mathsf   {COF}}}

\newcommand{\upannouncement}[1]{[\S\ref{#1} above]}

\newcommand{\cI}{\mathcal{I}}

\newcommand{\M}{\mathcal{M}}
\newcommand{\Null}{\mathcal{N}}

\newcommand{\add}{\mathsf{add}}

\newcommand{\R}{\mathbb{R}}
\newcommand{\EdNote}[1]{\par\medskip\noindent\textbf{#1.}}

\renewcommand{\b}{\mathfrak{b}}

\newcommand{\bq}{\begin{quote}}
\newcommand{\eq}{\end{quote}}
\renewcommand{\O}{\mathcal{O}}
\newcommand{\B}{\mathcal{B}}
\newcommand{\BG}{\B_\Gamma}

\newcommand{\BO}{\B_\Omega}

\newcommand{\sone}{\mathsf{S}_1}    \newcommand{\sfin}{\mathsf{S}_{fin}}
\newcommand{\ufin}{\mathsf{U}_{fin}}

\newcommand{\Union}{\bigcup}
\newcommand{\nin}{\not\in}
\newcommand{\cF}{\mathcal{F}}

\newcommand{\cU}{\mathcal{U}}

\newcommand{\fU}{\mathfrak{U}}

\newcommand{\NN}{{{}^{\naturals}\naturals}}
\newcommand{\naturals}{{\mathbb N}}
\newcommand{\N}{\naturals}

\newcommand{\sbst}{\subseteq}
\newcommand{\by}[2]{\par\hfill\emph{#1}, #2}

\newcommand{\Tau}{\mathrm{T}}
\newcommand{\CE}{\textsc{CE}}

\newcommand{\be}{\begin{enumerate}}
\newcommand{\ee}{\end{enumerate}}
\newcommand{\bi}{\begin{itemize}}
\newcommand{\ei}{\end{itemize}}
\renewcommand{\i}{\item}

\newcommand{\general}{\small\vfill\par\noindent\hrulefill\par
\noindent\textbf{Previous issues.}
The first issues of this bulletin,
which contain general information (first issue),
basic definitions, research announcements, and open problems (all issues) are available online,
on \arx{math.GN/$x$}, where $x$ is \texttt{0301011}, \texttt{0302062}, \texttt{0303057},
\texttt{0304087}, \texttt{0305367}, \texttt{0312140}, and \texttt{0401155},
respectively, for issues number $1$ to $7$.\\[0.1cm]
\textbf{Contributions.}
Please submit your contributions (announcements, discussions, and open problems)
by e-mailing us. It is preferred to write them
in \LaTeX{}.
The authors are urged to use as standard notation as possible, or otherwise give
the definitions or a reference to where the notation is explained.
Contributions to this bulletin would not require any transfer of copyright,
and material presented here can be published elsewhere.\\[0.1cm]
\textbf{Subscription.}
To receive this bulletin (free) to your
e-mailbox, e-mail us:\\
{tsaban@math.huji.ac.il}
}

\newcommand{\ArxPaper}[5]{\subsection{#3}{#5}\par\hfill{\arx{#1}}\par\hfill\emph{#4}\par\hfill{#2}}
\newcommand{\SPMBul}{\textbf{$\mathcal{SPM}$ Bulletin}}

\newcommand{\arx}[1]{\texttt{http://arxiv.org/abs/#1}}

\newcommand{\probmonth}{\emph{Problem of the month}}

\title[$\mathcal{SPM}$ Bulletin \textbf{8} (March 2004)]{%
$\mathcal{SPM}$ Bulletin\\[0.5cm]
Issue number 8: March 2004 \CE{}}

\begin{document}
\maketitle

\tableofcontents

\section{Editor's note}

Never has an issue of the \SPMBul{} contained as much interesting information
as this issue does.
In addition to the interesting research announcements, this issue contains
announcements of solutions for \emph{three}
open problems, one of which being a \probmonth{} in an earlier issue.
\medskip

Contributions to the next issue are, as always, welcome.

\medskip

\by{Boaz Tsaban}{tsaban@math.huji.ac.il}

\hfill \texttt{http://www.cs.biu.ac.il/\~{}tsaban}

\section{Research announcements}

\ArxPaper{math.GN/0312477}{lkocinac@ptt.yu}
{The Reznichenko property and the Pytkeev property in hyperspaces}
{Ljubisa D.\ R.\ Kocinac}
{We investigate two closure-type properties,
the Reznichenko property and the Pytkeev property,
in hyperspace topologies.}

\ArxPaper{math.LO/0401134}{zapletal@math.ufl.edu}
{Between Maharam's and von Neumann's problems}
{Ilijas Farah and Jindrich Zapletal}
{If $I$ is a suitably definable $\sigma$-ideal on the real
line and the factor algebra of Borel sets modulo $I$ is weakly
distributive then the algebra carries a Maharam submeasure.}

\ArxPaper{math.LO/0305241}{}
{Fraisse Limits, Ramsey Theory, and Topological Dynamics of Automorphism Groups}
{A.\ S.\ Kechris, V.\ G.\ Pestov and S.\ Todorcevic}
{We study in this paper some connections between the Fraisse
theory of amalgamation classes and ultrahomogeneous structures,
Ramsey theory, and topological dynamics of automorphism groups of countable structures.}

\ArxPaper{math.LO/0401343}
{don.monk@colorado.edu}
{Concerning problems about cardinal invariants on Boolean algebras}
{James Donald Monk}
{The present status of the problems in my book \emph{Cardinal Invariants on Boolean
algebras} (Birkhauser 1996) is described, with a description of solutions or
partial solutions, and references.}

\ArxPaper{math.GN/0402066}{}
{Small Locally Compact Linearly Lindelof Spaces}
{Kenneth Kunen}
{There is a locally compact Hausdorff space of weight $\aleph_\omega$ which is
linearly Lindelof and not Lindelof. This improves an earlier result, which
produced such a space of weight $\beth_\omega$.}

\subsection{Luzin gaps}
We isolate a class of $F_{\sigma\delta}$ ideals on $\mathbb{N}$  that includes
all analytic $P$-ideals  and all $F_\sigma$ ideals,  and introduce `Luzin gaps'
in their quotients. A dichotomy for Luzin gaps allows us to freeze gaps,  and
prove some gap preservation results.  Most importantly, under PFA all
isomorphisms between quotient algebras over these ideals have continuous
liftings. This gives a partial confirmation  to  the author's rigidity
conjecture for quotients $\mathcal{P}(\mathbb{N})/\mathcal{I}$.  We also prove
that the
ideals $\operatorname{NWD}(\mathbb{Q})$ and $\operatorname{NULL}(\mathbb{Q})$
have the Radon--Nikod\'ym property,  and (using OCA$_\infty$) a uniformization
result for $\mathcal{K}$-coherent families of continuous  partial functions.

\medskip
To appear in: \emph{Transactions of the AMS}.
\by{Ilijas Farah}{ifarah@mathstat.yorku.ca}

\subsection{The weak Fr\'echet-Urysohn property in function spaces}\label{sakaiwfu}
For a Tychonoff space $X$, we denote by $C_p(X)$
the space of all real-valued continuous functions on $X$ with
the topology of pointwise convergence.
In this paper, we note that for every analytic space $X$,
$C_p(X)$ is weakly Fr\'echet-Urysohn, and
solve a related problem of Tsaban \cite{reznicb}.
\by{Masami Sakai}{sakaim01@kanagawa-u.ac.jp}

\EdNote{The solved problem}
Following \cite{coc7}, we say that a cover $\cU$ of $X$ is \emph{$\fU$-groupable} if
there is a partition of $\cU$ into finite sets, $\cU = \Union_{n\in\N}\cF_n$, such  that
for each infinite subset $A$ of $\N$, $\{\cup\cF_n\}_{n\in A}\in\fU$.
Let $\fU^\gpbl$ be the family of $\fU$-groupable elements of $\fU$.

In  \cite{coc7} it is proved that the Hurewicz covering property $\ufin(\O,\Gamma)$
is equivalent to $\sfin(\Lambda,\Lambda^{gp})$.
In \cite{hureslaloms} it is proved that
$\sfin(\Lambda,\Lambda^{gp})=\binom{\Lambda}{\Lambda^{gp}}$.
In \cite{reznicb} we asked whether the analogue result for $\binom{\Omega}{\Omega^{gp}}$ is true,
namely, whether $\sfin(\Omega,\Omega^{gp})=\binom{\Omega}{\Omega^{gp}}$.

Sakai \upannouncement{sakaiwfu} gave a negative answer in the following strong sense:
He showed that the Baire space $\NN$ (and, consequently, any analytic space)
satisfies the stronger property $\binom{\BO}{\BO^{gp}}$. Recall that the Baire
space does not even satisfy $\sfin(\O,\O)$ (Menger's property), which is the
weakest property in the Scheepers diagram.
\by{Boaz Tsaban}{tsaban@math.huji.ac.il}

\subsection{$\cF$-Hurewicz spaces}
We investigate a generalization of spaces that satisfy the Hurewicz
covering property. In particular we are interested in characterization of
such spaces in terms of some properties of function spaces.
\by{F.\ Cammaroto}{camfil@unime.it}
\by{Lj.\ D.\ R.\ Ko\v{c}inac}{lkocinac@ptt.yu}
\by{G.\ Nordo}{nordo@dipmat.unime.it}

\subsection{Open problems in topology}
A new survey of the book \emph{Open Problems in Topology} appeared
in \emph{Topology and its Applications} \textbf{136} (2004)
37--85. The survey and the original book are available from the
Elsevier Mathematics portal
\bq
\texttt{http://www.elseviermathematics.com/}
\eq
(select \emph{Books} and go to the bottom of the page).

I have also completed editing
\emph{Problems from Topology Proceedings}. The book is available online at
\begin{quote}
\texttt{http://at.yorku.ca/i/a/a/z/10.htm}\\
\texttt{http://arXiv.org/abs/math.GN/0312456}
\end{quote}
This book consists of material originally appearing in the Problem
Section of the journal \emph{Topology Proceedings} since 1976 as
well as some other well-known problem lists in general topology
from the 1970's that have some connection to the journal. The
problems have been updated with current information on solutions
with bibliographic references. In particular, the book features
these collections:
\begin{itemize}
\item All contributed problems to the Problem Section, classified by
subject;
\item Eight Classic Problems by Peter J. Nyikos, including
information  from the two recent articles Twenty-five years later;
\item New Classic Problems from 1990;
\item Problems from Mary Ellen Rudin's Lecture notes in set-theoretic
topology (1975/7);
\item Problems from A.V. Arhangelskii's Structure and classification
of topological spaces and cardinal invariants (1978);
\item Continuum theory problems by Wayne Lewis (1983);
\item Problems in continuum theory by Janusz R. Prajs, including
essays by Charles L. Hagopian and Janusz J. Charatonik;
\item Classification of homogeneous continua by James T. Rogers, Jr.,
including material from survey articles of 1983 and 1989.
\end{itemize}

\by{Elliott Pearl}{elliott@at.yorku.ca}

{
\newcommand{\kl}{(\kappa,\lambda)}
\newcommand{\ghw}{\mathcal{G}^{\kappa}_{\lambda-1}}
\newcommand{\klm}{(\kappa,\lambda,<\mu)}
\newcommand{\bB}{\mathbb{B}}
\subsection{The hyper-weak distributive law and a related game in Boolean algebras}
The hyper-weak $\kl$-distributive law, formulated by Prikry, is a very weak
generalization of the $\klm$-distributive law.
We  define a related infinitary, two-player  game, called $\ghw$, and  give
connections between the hyper-weak $\kl$-distributive law and the existence of
winning strategies  for  the two players of $\ghw$, obtaining a game-theoretic
characterization of the hyper-weak $\kl$-distributive law  for many  pairs of
cardinals $\kappa,\lambda$, under GCH.  We then construct $\kappa^+$-Suslin
algebras for every infinite cardinal $\kappa$  in which, for each infinite
cardinal $\lambda\le\kappa$, neither player has a winning strategy for
$\mathcal{G}_{\lambda-1}^{\kappa}$.  This shows that the gap between the
strengths of the properties ``II wins $\mathcal{G}^{\kappa}_1(\infty)$ in $\bB$''
and ``the  $(\kappa,\infty)$-distributive law holds in $\bB$'' is consistently even
larger than was previously known.

\medskip

For a related work see: Natasha Dobrinen, \emph{Games and general distributive laws in Boolean algebras},
Proc.\ AMS \textbf{131} (2003), 309--318. (See also:
Natasha Dobrinen, \emph{Errata to `Games and general distributive laws in Boolean algebras'},
Proc.\ AMS \textbf{131} (2003), 2967--2968.)

\medskip

\by{James Cummings}{jcumming@andrew.cmu.edu}
\by{Natasha Dobrinen}{dobrinen@math.psu.edu}
}

\subsection{Problem of Issue 7 solved}\label{solved}
In the seventh issue we have posed the following conjecture.

\begin{conj}
If $X$ has strong measure zero and $|X|<\b$,
then all finite powers of $X$ have strong measure zero.
\end{conj}

Bartoszynski has found a combinatorial proof of the conjecture
when $X$ is a subsets of the Cantor space (to appear in \cite{prods}).
We have very recently found out that Scheepers has proved the
following more general result in \cite{smzpow}, for arbitrary metric spaces:
If $X$ has the Hurewicz property $\ufin(\O,\Gamma)$ and
strong measure zero, and $Y$ has strong measure zero,
then $X\x Y$ has strong measure zero.
A similar result follows from \cite{NSW}, as shown in \cite{prods}.

The following, though, remains open.

\begin{prob}
Assume that $X$ satisfies the Hurewicz property $\ufin(\O,\Gamma)$ and has strong measure zero
(the last property can be replaced by $\sone(\O,\O)$ or meager-additive).
Does it follow that all finite powers of $X$ satisfy $\ufin(\O,\Gamma)$?
\end{prob}

A positive answer would imply that the Gerlitz-Nagy property $(*)_{GN}=\sone(\Omega,\Lambda^{gp})$ is
preserved under taking finite powers, and that
$(*)_{GN}=\sone(\Omega,\Omega^{gp})$, see the definitions in \upannouncement{sakaiwfu}.

\by{Boaz Tsaban}{tsaban@math.huji.ac.il}

\subsection{Errata to: \emph{Hereditary topological diagonalizations and the Menger-Hurewicz Conjectures}}
Taras Banakh and Lubomyr Zdomsky have found a gap in our mentioned paper.
A revision of the paper has been made,
and the results promised in the abstract of the original paper are still proved.
Currently the paper is being thoroughly checked by a colleague.
Once verified for correctness, it will be re-posted to the Mathematics ArXiv.

\by{Tomek Bartoszynski}{tomek@diamond.boisestate.edu}
\by{Boaz Tsaban}{tsaban@math.huji.ac.il}

\subsection{The additivity number of the Menger and Scheepers properties}
For a family $\cI$ of subspaces of a topological space,
define $\add(\cI)=\min\{|\cF| : \cF\sbst\cI\mbox{ and }\cup\cF\nin\cI\}$.
I have recently obtained the following results.
For subspaces of a hereditarily Lindel\"of topological space:
\be
\i $\add(\ufin(\Gamma,\O))\ge\g$,
\i Under NCF, $\add(\ufin(\Gamma,\Omega))=\d$.
\ee
The first result gives a negative answer to Problem 2.4 from
\cite[full version]{huremen2}, which asks whether $\add(\ufin(\Gamma,\O))=\b$.
Indeed, under $\mathfrak u<\mathfrak g$ we have $\mathfrak b<\mathfrak g$.

The second result is a straightforward translation of some other results of
Taras Banakh, and it strengthens the result of \cite[full version]{huremen2},
that under NCF, $\add(\ufin(\Gamma,\Omega))\ge\max\{\b,\g\}$.

\by{Lubomyr Zdomsky}{lubomyr@opari.ltg.lviv.ua}

\section{Conferences}

\subsection{Geometric Topology: Infinite-Dimensional Topology, Absolute Extensors, Applications}
The conference will be held
on  26--29 of May 2004, at the Lviv Ivan Franko National University, Lviv, Ukraine.

Organizing/Program Committee:
T.\ Banakh, R.\ Cauty, A.\ Chigogidze, J.E.\ Keesling, V.\ Kyrylych, Ya.\ Prytula, D.\ Repovs, O.\ Skaskiv, M.\ Zarichnyi.

We plan plenary talks (45 min), section talks (25 min) and short communications
(10 min). The Conference will be held in the Main Building of the University in the historical
center of Lviv.
The official language of the Conference is English.
Persons interested in participating at the Conference are kindly asked to register by e-mail as
well as to send one page abstract in English, which should be prepared in LaTeX by e-mail till March 30, 2004.
The registration fee is 50 USD (for accompanying persons 25 USD), for participants from the countries of
the former Soviet Union the registration fee is 20 USD (for accompanying persons 10 D).
The registration fee is to be paid upon arrival in Lviv. The fee covers organization costs: abstracts,
tea/coffee during breaks, and cultural program.

The accomodation price in Lviv at the moment is from 10 to 80 USD per night
depending on facilities.
For further information, please visit the following web-sites:
\bq
\verb|http://www.all-hotels.com.ua/addz2.php3?Lang=1&sCity=4|\\
\verb|http://www.piligrim.lviv.ua/ukraine/page6_en.html|
\eq
The Organizing Committee will have an opportunity to pre-book rooms in Students Hostel
(for interested participants).
We regret that travel and daily expenses cannot be paid by the Organizing Committee.

For further information or specific requests please e-mail us to the following address.

\by{Oleg V.\ Gutik}{ogutik@iapmm.lviv.ua}

\section{Problem of the month}

For an ideal ${\mathcal J}\sbst P(\R)$, we say that
$H\sbst\R\x\R$ is a \emph{Borel ${\mathcal J}$-set} if
$(H)_x \in {\mathcal J}$ for all $x \in \R.$

Using this terminology we define the following classes of small sets:

\begin{itemize}
\item[$\COF(\mathcal{J})$] $=\{X \subseteq \R: \mbox{for every Borel $\mathcal{J}$-set }H, \{(H)_x: x \in X\}
\mbox{ is not a basis of }\mathcal{J}\}$,
\item[$\ADD(\mathcal{J})$] $=\{X \subseteq \R: \mbox{for every Borel $\mathcal{J}$-set  }H,
\bigcup_{x \in X} (H)_x \in \mathcal{J}\}$,
\item[$\COV(\mathcal{J})$] $=\{X \subseteq \R: \mbox{for every Borel $\mathcal{J}$-set }H,
\bigcup_{x \in X} (H)_x \neq \R\}$,
\item[$\NON(\mathcal{J})$] $=\{X \subseteq \R:\mbox{every image of $X$ by a Borel function is in }\mathcal{J} \}$.
\item[${\mathsf D}$] $=\{ X\sbst \R : \mbox{ for every Borel function }f:\R\to\NN\mbox{, }f[X]\mbox{ is not
dominating}\}$.
\item[${\mathsf B}$] $=\{ X\sbst \R : \mbox{ for every Borel function }f:\R\to\NN\mbox{, }f[X]\mbox{ is bounded}\}$.
\end{itemize}

\begin{figure}[!h]
\unitlength=.95mm
\begin{picture}(140.00,60.00)(10,10)
\put(20.00,20.00){\makebox(0,0)[cc]
{${\sf ADD}(\Null)$ }}
\put(60.00,20.00){\makebox(0,0)[cc]
{${\sf ADD}(\M)$ }}
\put(100.00,20.00){\makebox(0,0)[cc]
{${\sf COV}(\M)$ }}
\put(140.00,20.00){\makebox(0,0)[cc]
{${\sf NON}(\Null)$ }}
\put(20.00,60.00){\makebox(0,0)[cc]
{${\sf COV}(\Null)$ }}
\put(60.00,60.00){\makebox(0,0)[cc]
{${\sf NON}(\M)$ }}
\put(100.00,60.00){\makebox(0,0)[cc]
{${\sf COF}(\M)$ }}
\put(140.00,60.00){\makebox(0,0)[cc]
{${\sf COF}(\Null)$ }}
\put(100.00,40.00){\makebox(0,0)[cc]
{${\sf D}$ }}
\put(60.00,40.00){\makebox(0,0)[cc]
{${\sf B}$ }}
\put(28.00,61.00){\vector(1,0){20.00}}
\put(110.00,20.00){\vector(1,0){20.00}}
\put(110.00,61.00){\vector(1,0){20.00}}
\put(70.00,61.00){\vector(1,0){20.00}}
\put(28.00,20.00){\vector(1,0){20.00}}
\put(70.00,20.00){\vector(1,0){20.00}}
\put(70.00,40.00){\vector(1,0){20.00}}

\put(20.00,25.00){\vector(0,1){29.00}}
\put(60.00,25.00){\vector(0,1){11.00}}
\put(60.00,42.00){\vector(0,1){12.00}}
\put(100.00,25.00){\vector(0,1){11.00}}
\put(100.00,42.00){\vector(0,1){12.00}}
\put(140.00,25.00){\vector(0,1){29.00}}
\end{picture}
\caption{Chichon Diagram for small sets of reals \label{chichondiagr}}
\end{figure}
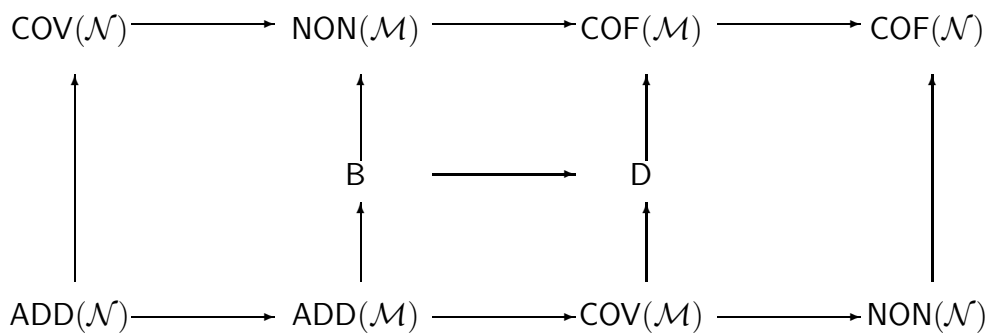

The interrelationships of these classes were extensively studied in \cite{pawlikowskireclaw}
in the case that $\mathcal{J}$ is $\M$ (the ideal of meager sets of reals)
or $\Null$ (the ideal of Lebesgue measure zero sets of reals).
These are summarized in Figure \ref{chichondiagr}.

The relationship between Figure \ref{chichondiagr} and the well-known Chichon diagram that expresses provable
relationships among certain cardinal numbers is that a cardinal number in a particular position in Cichon's
diagram is actually the minimal cardinality for a set of real numbers not belonging to the class in the
corresponding position in Figure \ref{chichondiagr}.

Many of these classes can be defined by selection principles, e.g.:
$\mathsf{D}=\sone(\BG,\B)$, $\mathsf{B}=\sone(\BG,\BG)$,
$\COV(\M)=\sone(\B,\B)$,
$\ADD(\M)=\sone(\B,\B)\cap\sone(\BG,\BG)$.

\medskip

The following problem was suggested by Bartoszynski:
It is known \cite{pawlikowskireclaw}
that if $X \nin \NON(\M)$ and $Y\nin\sone(\BG,\B)$, then $X\x Y\nin \COF(\M)$.

\begin{prob}
Suppose that $X \nin \NON(\M)$ and $Y\nin\mathsf{D}$.
Does it imply that $X\cup Y\nin \COF(\M)$?
\end{prob}

\by{Tomek Bartoszynski}{tomek@diamond.boisestate.edu}

\section{Problems from earlier issues}
In this section we list the past problems posed in the \SPMBul{},
in the section \probmonth{}.
For definitions, motivation and related results, consult the
corresponding issue.

For conciseness, we make the convention that
all spaces in question are
zero-dimentional, separable metrizble spaces.

\begin{issue}
Is $\binom{\Omega}{\Gamma}=\binom{\Omega}{\Tau}$?
\end{issue}

\begin{issue}
Is $\ufin(\Gamma,\Omega)=\sfin(\Gamma,\Omega)$?
And if not, does $\ufin(\Gamma,\Gamma)$ imply
$\sfin(\Gamma,\Omega)$?
\end{issue}

\begin{issue}
Does there exist (in ZFC) a set satisfying
$\ufin(\O,\O)$ but not $\ufin(\O,\Gamma)$?
\end{issue}
\begin{proof}[Solution]
\textbf{Yes} (Lubomyr Zdomsky, 2003).
\end{proof}

\begin{issue}
Does $\sone(\Omega,\Tau)$ imply $\ufin(\Gamma,\Gamma)$?
\end{issue}

\begin{issue}
Is $\p=\p^*$?
\end{issue}

\begin{issue}
Does there exist (in ZFC) an uncountable set satisfying $\sone(\BG,\B)$?
\end{issue}

\begin{issue}
Assume that $X$ has strong measure zero and $|X|<\b$.
Must all finite powers of $X$ have strong measure zero?
\end{issue}
\begin{proof}[Solution]
\textbf{Yes} (Scheepers \cite{smzpow}; Bartoszynski). See \upannouncement{solved}.
\end{proof}


\general


\begin{thebibliography}{00}

\bibitem{huremen2}
T.\ Bartoszynski, S.\ Shelah, and B.\ Tsaban,
\emph{Additivity properties of topological diagonalizations},
The Journal of Symbolic Logic \textbf{68} (2003),
1254--1260.
(Full version: \arx{math.LO/0112262})

\bibitem{coc7}
Lj.\ D.\ R.\ Ko\v{c}inac and M.\ Scheepers,
\emph{Combinatorics of open covers (VII): groupability},
Fundamenta Mathematicae \textbf{179} (2003),
131--155.

\bibitem{NSW}
A.\ Nowik, M.\ Scheepers, and T.\ Weiss,
\emph{The algebraic sum of sets of real numbers with strong measure zero sets},
J.\ Symbolic Logic \textbf{63} (1998), 301--324.

\bibitem{pawlikowskireclaw}
J.\ Pawlikowski and I.\ Rec{\l}aw,
\emph{Parametrized Cicho\'n's diagram and small sets},
Fundamenta Mathematicae \textbf{147} (1995),
135--155.

\bibitem{smzpow}
M.\ Scheepers,
\emph{Finite powers of strong measure zero sets},
The Journal of Symbolic Logic \textbf{64} (1999),
1295--1306.

\bibitem{reznicb}
B.\ Tsaban,
\emph{The minimal cardinality where the Reznichenko property fails},
Israel Journal of Mathematics \textbf{140} (2004),
367--374.
\arx{math.GN/0304024}

\bibitem{hureslaloms}
B.\ Tsaban,
\emph{The Hurewicz covering property and slaloms in the Baire space},
submitted.
\arx{math.GN/0301085}

\bibitem{prods}
B.\ Tsaban and T.\ Weiss, \emph{Products of special sets of real
numbers},
eprint \arx{math.LO/0307226}

\end{thebibliography}
\end{document}